\documentclass[journal]{IEEEtran}
\usepackage[utf8]{inputenc}
\usepackage{graphicx}
\usepackage{cite,url}
\usepackage[T1]{fontenc}
\usepackage{graphicx}
\usepackage{amssymb,amsmath,amstext}
\usepackage{epstopdf}
\usepackage{listings,color}
\usepackage{fp,float}
\usepackage{tikz,tikz-inet,pgf,pgfplots}
\usepackage{algorithm}
\usepackage[noend]{algpseudocode}
\usepackage{etoolbox}

\usepackage{ifthen}

\usetikzlibrary{patterns}

\makeatletter
\def\BState{\State\hskip-\ALG@thistlm}
\makeatother

\newtheorem{assumption}{Assumption}

\newtheorem{proposition}{Proposition}

\definecolor{colorOn}{RGB}{200,200,200}

\newbool{short}
\booltrue{short}

\newcommand\short[2]{\ifbool{short}{#1}{#2}}

\newcommand\hpat[1]{}
\newcommand\define{:=}
\gdef\dd{{\mathrm d}}

\gdef\extravs{}

\pgfcreateplotcyclelist{\mylist}{
	{blue,mark=none,thick,smooth},
	{black,mark=none,thick,densely dashed,smooth},
	{red,mark=none,thick,densely dotted},
	{dash pattern=on 8pt off 2 pt on 2pt off 2pt,thick,mark=none},
}

\color{black}

\title{Decentralized demand response for temperature-constrained appliances}

\author{
    \IEEEauthorblockN{Nguyen Hoang Son Duong\IEEEauthorrefmark{1}, Patrick Maill\'e\IEEEauthorrefmark{1}, Ashish Kumar Ram\IEEEauthorrefmark{3}, Laurent Toutain\IEEEauthorrefmark{1}}\\
    \IEEEauthorblockA{%
 \begin{tabular}{ccc}
 \IEEEauthorrefmark{1}IMT Atlantique&~~~~~~~~~~~&\IEEEauthorrefmark{2}Technische Universit\"at Berlin\\
Rennes, FRANCE&&Berlin, GERMANY\\
\{first.last\}@imt.fr&&ashish.ram88@gmail.com
\end{tabular}
}
}

\begin{document}
\maketitle

\begin{abstract}
The evolution of the power grid towards the so-called Smart Grid, where information technologies help improve the efficiency of electricity production, distribution and consumption, allows to use the fine-grained control brought by the Internet of Things capabilities 
to perform distributed demand response when requested by the grid.

In this paper, we study the demand-response potential of coordinated large numbers of appliances which have to maintain some temperature within a fixed range through the ON/OFF functioning of a temperature modifier. 
We introduce a mathematical model and methods to coordinate appliances with given requirements, in order to offer a global energy demand reduction for a desirable duration while still satisfying the temperature constraints, and with limited communication overhead. We quantify the maximum power reduction that can be attained, as a function of the reduction duration asked by the grid.
\end{abstract}

\begin{IEEEkeywords}
	Smart Grid, Demand Response (DR), temperature constraints, Internet of Things, communication technologies\end{IEEEkeywords}

\IEEEpeerreviewmaketitle

\section{Introduction} \label{sec:introduction}

The transition toward the Smart Grid paradigm is driven by several forces, including the raise and evolution of electricity demand (in particular, due to the appearance of electric vehicles~\cite{shuai2016charging}), the limited capacities of current grids~\cite{farina2015domino}, and environmental as well as economic considerations leading to a strong development of renewable energy sources. 
Power systems will largely rely on information and communication technologies to optimize and coordinate the production, transmission, distribution, and consumption of electricity, to improve efficiency and reliability~\cite{farhangi2010path,gungor2011smart}. 

Among the new important aspects of the Smart Grid is the increasing need for flexibility at the consumption side: indeed, with the new constraints imposed by renewable energy production (intermittence, uncontrollability, and only partial predictability) and the difficulties for storing electricity (costs and limited efficiency), an interesting direction to maintain the balance between production and consumption is to affect demand based on the grid conditions. This is the so-called demand response (DR) approach~\cite{rahimi2010demand,albadi2007demand}, which can consist in shifting demand in time, or in rewarding users who accept to adapt their consumption when asked to. 
Through DR, the grid partially controls the consumption of electric appliances. This can be done with an EMS  (Energy Management Server) using data delivered through the AMI (Advanced Metering Infrastructure, involving networking technologies) from and to the appliances~\cite{lee2010design,logenthiran2012demand,mohsenian2010autonomous}. 
Demand reduction requests are issued by utilities or by market managers (which will be indifferently called ``the grid'' thereafter), either directly to customers, or more likely to \emph{aggregators}, which are new entities in the electricity market behaving like brokers between several users and the utility operator. Aggregators contract with several consumers with flexibility potential, and coordinate them to offer significant-scale flexibility offers
~\cite{gkatzikis2013role}.

Although demand response is already applied in the grid, it still concerns only large consumers, with whom it is simpler to establish contracts for flexibility services (typically, demand reduction during some peak periods). But current research aims at leveraging the demand response potential of smaller consumers like individual households, which raises several computational issues~\cite{gilbert2015optimal}. 
In this paper, we also apply DR with many small consumers, but for a specific set of electricity-consuming appliances such as fridges, A/C systems, or water heaters. Such appliances exist in very large numbers, so that if coordinated they can offer significant demand-response services to the grid.
The most difficult challenge with these appliances is that we cannot turn them off for an undetermined duration, because they are subject to temperature constraints: the temperature has to remain within a predefined range. In this paper, we present a simple mathematical model and its analysis to coordinate such appliances in order to offer a global energy demand reduction for a given duration while still satisfying the temperature constraints. In practice, the grid would request a power reduction for a given duration, and an aggregator controlling a large number of those appliances could respond by coordinating these appliances using our schemes. 

\short{}{Note also that the flexibility market is associated with incentive prices, since the actors have to be sufficiently incentivized to provide a flexibility service. That aspect is out of the scope of this paper: we rather focus on the amount of flexibility that an aggregator can offer given the number of appliances it controls. That specific number would depend on the contract prices, hence our results could be applied to an economic model in order to estimate the appropriate amounts. Such considerations are left for future work, since they involve many additional assumptions to be made, on the flexibility market, on the structure of the aggregator market, and on user preferences.}

Our main results are based on mathematical models for the appliance behaviors, and their analysis in terms of reduced power or saved energy. The advantage of our proposition is that the communication overhead remains very small (i.e., only one broadcast message from the aggregator to the appliances, upon the aggregator receiving a request from the grid), a desirable feature in IoT. Furthermore, we provide several methods for implementing a power reduction, depending on the grid needs. \emph{(i)} If the duration and the amplitude of the asked power reduction are limited, then a very basic mechanism is sufficient to satisfy the request. \emph{(ii)} a slightly more elaborate mechanism can be implemented for longer and/or bigger requests. Possibly, there is still some room for further optimization but \emph{(iii)} we compute upper bounds for the flexibility service that can be satisfied with those appliances, and show how our two simple proposals perform with respect to that upper bound. 

Finally, while all our results are derived for consumption reduction requests, they can be transposed to providing consumption increases, a service for which requests are more rare, but occur in practice~\cite{gotz2013negative} and may occur more frequently due to the increase of off-peak production from renewable energy sources (e.g., from wind farms at night). 

The remainder of this paper is organized as follows. Section \ref{sec:model} describes the mathematical model for the appliance behavior, and the format of reduction requests issued by the grid. The upper bound for the relative power reduction that can be offered over a given duration is computed in Section \ref{sec:upper_bound}, where we ignore constraints on the shape of the reduction over time (the reduction is not constant over this duration). To include the constraint of constant power reduction, a simple mechanism is presented in Section \ref{sec:schOne}, which is further improved in Section \ref{sec:schTwo} to get closer to the upper bound. Section \ref{sec:conclusion_perspectives} concludes the paper, suggesting more research and perspectives.  

\section{Mathematical Model}\label{sec:model}

\subsection{Appliance behavior}

We consider appliances that have to maintain some temperature within a fixed range $[T_{\min},T_{\max}]$ of size $\Delta\define T_{\max}-T_{\min}$, through the ON/OFF functioning of a temperature modifier which consumes some power $P$ when ON and no power otherwise. 
Hence our model applies to heating or AC systems, fridges/freezers, water heaters, etc.
In the figures displaying temperatures, cooling appliances are considered, but our model is generic, and all our mathematical formulations are agnostic to the type of appliance (heating or cooling).

When the temperature modifier is ON, we assume the temperature varies (increases for a heating appliance and decreases for a cooling appliance) with constant speed $v$ (degrees per time unit); otherwise it drifts in the opposite direction with constant speed $w$.
Finally, that temperature modifier is only turned ON when necessary, i.e., for a heating (resp., cooling) appliance, when the temperature has drifted to the lower (resp., upper) limit of the interval $[T_{\min},T_{\max}]$. It then remains on use until the upper (resp., lower) limit of that interval is reached.

In this paper, we assume that all appliances considered \short{}{are identical, in the sense that they all }have the exact same characteristics. This is actually without loss of generality, since we can deal with any finite number of appliance categories (including, mixing cooling and heating appliances) by simply aggregating the flexibility possibilities of all categories. The aim of this paper is to compute the flexibility potential of one given category, when the appliances in that category follow the behavior  described above.

Summarizing, we have the following assumption:\extravs
\begin{assumption}[Individual nominal appliance behavior] \label{assumption:Delta_v_w}
Without reduction requests, each appliance functions through \emph{cycles} of total duration $\Delta (1/v + 1/w)$, during which the temperature modifier is ON and consumes some power $P$ for $\Delta/v$ time units, and consumes no energy for $\Delta/w$.
\end{assumption}
\extravs

Note that for simplicity, we ignore here some possible extra consumption costs upon launching the engine; incorporating this into the model, as well as more complex consumption patterns over a cycle, is left for future work.

\subsection{Desynchronization among appliances}

We make the reasonable assumption that appliances are desynchronized, i.e., the points where they are in their cycles are uncorrelated. Assuming a large number of appliances, which we thereafter treat as a continuum, we end up with a uniform distribution of appliance positions (with respect to the cycle origin) over the cycle duration $\Delta (1/v + 1/w)$.

Note that the reduction requests issued by the grid will affect the appliances and destroy this uniform distribution. But in practice, reduction requests occur very rarely with respect to the cycle duration, so we can consider that thanks to uncorrelated small variations among appliance cycles (due to external causes such as user actions), a steady-state desynchronized situation with the uniform distribution is reached between two consecutive requests. This is summarized below.

\extravs
\begin{assumption} \label{assumption:demand_reduction_request}
When a demand reduction request is issued, all appliances are desynchronized: among appliances, the times $y$ since the beginning of their cycle are uniformly distributed over the interval $[0,\Delta (1/v+1/w)]$. 
\end{assumption}
\extravs

That steady-state situation is illustrated in Figure \ref{fig:analysis_simple}, showing the evolution of ON/OFF states of appliances over time. 


\begin{figure} [!h]
\def\pplus{\!\!+\!\!}
{\footnotesize
	\begin{tikzpicture}[scale=.32]

	\fill[colorOn]  (0,0) -- (0,4)--(4,0)--cycle;
	\fill[colorOn]  (0,9) -- (4,9)--(13,0)--(9,0)--cycle;
	\fill[colorOn]  (9,9) -- (13,9)--(22,0)--(18,0)--cycle;
	\fill[colorOn]  (18,9) -- (22,9)--(23,8)--(23,4)--cycle;

	\draw[thick,green] (0,4) node[left,black] {$\frac{\Delta}{v}$}--(4,0) node [below,black] {$\frac{\Delta}{v}$};
	
	\draw[thick,red] (0,9) node [left,black] {$\frac{\Delta}{v}\pplus\frac{\Delta}{w}$} --(9,0) node [below,black]{$\frac{\Delta}{v}\pplus\frac{\Delta}{w}$};	
	\draw[thick,red] (9,9) --(18,0) node [below,black]{$\frac{2\Delta}{v}\pplus\frac{2\Delta}{w}$};
	\draw[thick,red] (18,9) --(23,4);
	
	\draw[thick,green] (4,9)--(13,0);
	\draw[thick,green] (22,0) -- (13,9);
	\draw[thick,green] (22,9) -- (23,8);
	
	\node [font=\fontsize{8}{144}\selectfont]at (13,-1){$\frac{2\Delta}{v}\pplus\frac{\Delta}{w}$};
		
	\node[font=\fontsize{8}{144}\selectfont] at (1,1) {\emph{ON}};
	\node[font=\fontsize{8}{144}\selectfont] at (10,1) {\emph{ON}};
	\node[font=\fontsize{8}{144}\selectfont] at (19,1) {\emph{ON}};
	
	\node[red] at (5.5,1) {\emph{OFF}};
	\node[red] at (14.5,1) {\emph{OFF}};
	\node[red] at (23.5,1) {\emph{OFF}};
	
	\draw[->] (0,0) -- (23,0) node[below] {Time};
	\draw[->] (0,0) -- (0,11) node[left] {$y$};	
	\node[left] at (-0.25,0) {0};
	\draw[thick,densely dashed](0,9)--(23,9);

	\draw (-2,4.5)node[rotate=90]{Appliances};
	\end{tikzpicture}
	}
	\caption{The ON/OFF states of appliances over time. The behavior of a specific appliance should be read as a horizontal line on the graph, the height $y$ of that line (in the range $[0,\Delta/v+\Delta/w]$) being the time since the temperature-modifying system was last turned ON at time $0$. That value $y$ is specific to each appliance, and distributed uniformly among appliances.} 
	\label{fig:analysis_simple} 
\end{figure}

In the desynchronized steady-state, the proportion of appliances in ON state (and consuming $P$) always equals $\frac{1/v}{1/v+1/w}=\frac{w}{v+w}$. Hence, denoting by $N$ the number of appliances (assumed large), the aggregated consumption\short{}{ from those appliances} is constant and equals
\gdef\Ptot{P_{\text{tot}}}
\begin{equation}\label{eq:power_steadystate}
\Ptot = N P \frac{w}{v+w}.
\end{equation}

\subsection{Consumption reduction requests}

When the grid needs a reduction in the aggregated consumption, it sends a request to the aggregator, which will then coordinate the appliances spread over the territory. 
There are two main components in a reduction request, namely:
\begin{itemize}
\item the {\bf duration} over which the reduction should take place, that we will denote by $t$;
\item the {\bf amplitude} of the reduction, that is, the power reduction (in watts) over that duration.
\end{itemize}
We focus on deciding whether the aggregator can satisfy such a request. To do so, we compute for each possible duration, the maximum amplitude that can be offered by coordinating the appliances: if the result exceeds the asked amplitude, the aggregator will be able to satisfy the request\short{}{, providing the exact amplitude needed by having only a proportion of the appliances reduce consumption}.

Since the absolute amplitude value depends on the number of appliances responding to the aggregator, we rather reason in relative values. Similarly, we take as the reference value the consumption $\Ptot$ at the steady-state, expressed in~\eqref{eq:power_steadystate}, since this is the maximum possible reduction. 
Hence in the following, an amplitude reduction $R(t)$ over a period $t$ means that the average power reduction equals $R(t)\Ptot = R(t)\times N P \frac{w}{v+w}$.

Note that there can be other considerations in a reduction request. In particular, a demand reduction often consists in shifting demand in time, and therefore leads to a demand increase after the reduction period. This phenomenon is usually called the \emph{rebound effect}, and can exceed the reduction offered (this is not the case with our model). In this paper, we do not analyze what happens after the end of the reduction period: we assume that the consumption peak is over and that the grid can cope with the possible extra consumption resulting from the reduction. However, with our model that rebound effect is easy to compute numerically.

It is also natural to assume that the grid expects a \emph{constant} power reduction over the requested duration; we relax that assumption in the next section to compute an upper bound of the reduction amplitude, but include it in the other sections.

\section{Upper Bound for the consumption reduction} \label{sec:upper_bound}

In this section, we assume that the system aims at maximizing the aggregated energy (or equivalently, the average power) saved over the reduction duration $t$, ignoring any constraint on the consumed (or reduced) power to be constant over that interval. This will provide us with an upper bound for the reduction amplitude the aggregator can provide to the grid. 

\extravs
\begin{proposition} \label{prop:upper_bound}
Over a duration $t$, the relative reduced power with respect to the stationary consumption cannot exceed
\begin{equation}\label{eq:RedUpper_bound}
R_{\sup}(t) \define \max\left(1 - \frac{wt}{2\Delta}, \frac{\Delta}{2wt}\right)
\end{equation}
\end{proposition}
\extravs

\begin{IEEEproof}
Since we do not need the consumption to be constant over $t$ here, we can imagine a simple mechanism to minimize the total energy saved: just turn off the temperature-modifier upon receiving the request, let the temperature drift until the (upper for a cooling appliance, lower for a heating appliance) limit, and perform very short ON-OFF stages just to maintain that limit temperature\short{}{, as illustrated in Figure~\ref{fig:illustr_upper_bound}}. Given the parameters, in that second phase the proportion of time in ON stage should equal $\frac{w}{v+w}$.
\short{}{\begin{figure} [htbp]
\centering
{\footnotesize
\begin{tikzpicture}[xscale=1.4]
\draw[->] (0,0) -- (4,0) node[right] {Time};
\draw[->] (0,0) -- (0,4) node[above] {Temperature};
\draw[dashed] (4,3)--(0,3);
\node[left] at (0,3) {$T_{\max}$};
\draw[dashed] (4,1)--(0,1);
\node[left] at (0,1) {$T_{\min}$};
\node[below] at (0,0) {$0$};
\node[left] at (0,1.3){$T$};
\draw[ultra thick] (0,1.3)--node[below,rotate=33]{slope $w$}(2,3);
\draw[ultra thick] (2,3)--(3.7,3);
\node[] at (1,3.2) {$\delta_T / w$};
\draw[<->,thick](2.5,1.3)--node[left]{$\delta_T$}(2.5,3);
\draw[densely dotted](3.7,0)--(3.7,3);
\draw[densely dotted](0,1.3)--(3.7,1.3);
\node[below] at (3.7,0){$t$};
\end{tikzpicture}
}
\caption{Temperature evolution during reduction duration $t$ for a cooling appliance, minimizing the energy consumption over $t$. The time origin is set to the moment the request is received.} 
\label{fig:illustr_upper_bound} 
\end{figure}}
Hence the consumption for each appliance is $0$ until the temperature limit is reached, and $P\frac{w}{v+w}$ afterwards.

We assume the appliances are in the steady-state when the request arrives, so that their positions in their cycle are uniformly distributed over the interval $[0,\Delta/v+\Delta/w]$. As a consequence, the appliance temperature $T$ also follows a uniform distribution over $[T_{\min},T_{\max}]$, and the temperature distance  $\delta_T$ to the limit follows a uniform distribution over $[0,\Delta]$. (For a cooling appliance, $\delta_T = T_{\max}-T$ whereas $\delta_T = T-T_{\min}$ for a heating appliance.)

With our mechanism, the time before reaching the temperature limit is $\delta_T/w$: an appliance does not consume energy during that time, and then consumes a power $P\frac{w}{v+w}$ until the end of the reduction period. Hence a total consumed energy
$
P\frac{w}{v+w}\left[t-\delta_T/w\right]^+,$
where $[x]^+\define \max(0,x)$.

\gdef\Emin{{\mathbf E}_{\min}}
Using the uniform distribution for $\delta_t$, the expected energy $\Emin$ consumed per appliance over $t$ is
\begin{eqnarray*}
\Emin &=& {\mathbb E}_{\delta_T}\!\!\left[P\frac{w}{v+w}\left[t-\delta_T/w\right]^+\right]\\
&=&\frac{P}{v+w}\! \int_{x=0}^{\Delta}\frac{\left[wt-x\right]^+}{\Delta}\dd x.\\
\short{}{&=&\frac{P}{\Delta(v+w)}\times\left\{\begin{array}{ll}
\frac{(wt)^2}{2}&\text{if }wt<\Delta\\
\Delta(wt-\frac{\Delta}{2})&\text{otherwise}
\end{array}\right.\\}
&=&\frac{Pw}{v+w}\times\min\!\left(\frac{wt^2}{2\Delta}\ ,\ t-\frac{\Delta}{2w}\right)\\
&=&tP\frac{w}{v+w}\times\min\!\left(\frac{wt}{2\Delta}\ ,\ 1-\frac{\Delta}{2wt}\right)\short{.}{,}
\end{eqnarray*}
\short{}{which is illustrated in Figure~\ref{fig:Emin_vs_t} as a function of the reduction duration $t$.
\begin{figure}[htbp]
\centering
{\footnotesize
\begin{tikzpicture}[xscale=1.4,yscale=.7]
\draw[->] (0,0) -- (5,0) node[right] {$t$};
\draw[->] (0,0) -- (0,5) node[above] {Energy};
\draw[line width=4pt,domain=0:2,smooth,opacity=.3] plot (\x,{0.5*(\x)^2});
\draw[red,thick, densely dashed,domain=0:3.1,smooth] plot (\x,{0.5*(\x)^2}) node[left] {$\frac{P}{v+w} \frac{w^2t^2}{2\Delta}$};
\draw[line width=4pt,opacity=.3] (2,2) -- (3.5,4.7) node[below right]{$\Emin$};
\node[blue] at (1.1,-1) {$\frac{Pw}{v+w} \left(t-\frac{\Delta}{2w}\right)$};
\draw[blue,thick](3.5,4.7)--(0.2,-1.24); 
\draw[dashed] (2,2)--(0,2);
\draw[dashed] (2,2)--(2,0);
\node[below] at (2,0) {$\frac{\Delta}{w}$};
\node[left] at (0,2) {$\frac{\Delta^2}{2}$};
\node[below] at (0,0) {$0$};
\end{tikzpicture}
}
\caption{Minimum possible consumed energy per appliance $\Emin$ over a duration $t$ (in bold), for a request initiated in the steady-state situation.} 
\label{fig:Emin_vs_t} 
\end{figure}}
In the last expression for $\Emin$, the first term $tP\frac{w}{v+w}$ is the average consumed energy per appliance in the steady-state regime over a duration $t$, hence the relative reduction is directly one minus the second term\short{,}{. This yields a maximum relative reduction over $t$
\begin{eqnarray*}
R_{\sup}(t) \short{}{&=& 1 -  \min\left(\frac{wt}{2\Delta},1-\frac{\Delta}{2wt}\right)\\}
&=&\max\left(1-\frac{wt}{2\Delta},\frac{\Delta}{2wt}\right),
\end{eqnarray*}}
giving the proposition.
\end{IEEEproof}

\short{}{That upper bound for the relative reduction is illustrated in Figure~\ref{fig:power}, where it is also compared to our first proposed mechanism guaranteeing a constant power reduction over the duration $t$, introduced in the next section.}

\gdef\schOne{IndivRed}
\gdef\schTwo{CoordRed}
\section{\schOne: A simple mechanism with constant power reduction} \label{sec:schOne}

From now on, we look for a way to offer a \emph{constant} power reduction over the duration $t$ requested by the grid, a constraint ignored when computing the bound  in Proposition~\ref{prop:upper_bound}.

In this section, we take the most simple approach, when there is no real coordination among appliances: the aggregator just mobilizes each individual appliance that can offer a reduction over the duration $t$. We will call \emph{\schOne{}} the corresponding mechanism. To compute the possible (relative) reduction amplitude with such a mechanism, in what follows we express the proportion of such appliances.

\subsection{Relying on individual duration-$t$ reductions}

Given a duration $t$, we investigate here the conditions under which an individual appliance can offer a constant consumption reduction over $t$.  The reduction being \emph{with respect to the no-request situation}, we take that situation as our reference. 

Figure~\ref{fig:effacement_vs_cycle} illustrates what an appliance can offer, depending on its position in its cycle upon receiving the request (i.e., the time $y$ since its temperature modifier was last switched ON). As stated previously, $y$ is in the interval $[0,\Delta/v+\Delta/w]$, and the temperature modifier is only ON if $y<\Delta/v$.
Naturally, an appliance in OFF state cannot reduce its consumption since it is not currently consuming. To offer a constant reduction over $t$, an appliance must satisfy two conditions:
\par\noindent
\begin{enumerate}
\item without the request it would have been ON during $t$; 
\item it can afford to be OFF instead during $t$, without its temperature exiting the allowed range $[T_{\min},T_{\max}]$.
\end{enumerate}

\subsection{How much can we reduce during $t$ with \schOne{}?}

\begin{figure}[htbp] 
	\centering 
	{\footnotesize
	\begin{tikzpicture}[xscale=1.2,yscale=.5]
	\coordinate (a) at (0,0);
	\draw [->] (a) -- +(0, 5) node[above,font=\fontsize{8}{144}\selectfont]{Temperature}; 
	\draw [->] (a) -- +(5, 0) node[anchor=north,font=\fontsize{8}{144}\selectfont]{Time}; 
	\draw [dotted]  ([yshift=4cm] a) coordinate (b) -- +(5, 0); 
	\draw (b) node (Tmax) [left] {$T_{\max}$};
	\draw [dotted]  ([yshift=1cm] a) coordinate (c) -- +(5, 0); 
	\draw (c) node (Tmin) [left] {$T_{\min}$};
	\draw[<->] (b-|Tmax.west)--node[left]{$\Delta$}(c-|Tmax.west);
	
	\draw[ultra thick] ([yshift=4cm] a) -- node [rotate=-32, below,near end] {$-v$} ++(2, -3) coordinate (motorOff) -- node [rotate=28, above,near end] {$w$}++(3, 3); 
	
	\coordinate (d) at ([yshift=-2cm]a);
	\draw [->] (d) -- coordinate [very near end] (f)  +(0, 1) ; 
	\draw [->] (d) -- +(5, 0); 
	
	\draw [fill=gray]  (f) -- (f -| motorOff) coordinate (g) -- (g |- d) -- (d) -- cycle ; 
	\draw [loosely dashed] (motorOff) -- (motorOff |- a) node(e)[anchor=north]{$\Delta/v$}; 
	\draw [loosely dashed] (e) -- (motorOff |- d) node(e)[anchor=north]{}; 
	\draw (d)node[anchor=south east,font=\fontsize{8}{144}\selectfont]{Consumption, no red.};
	
	\color{red}{
		\draw ([yshift=4cm] a)++(0.22,-.33) node (aux1){};
		\draw[ultra thick,densely dashed] (aux1.center)--++(.33,.33)node(b1){}--++(2,-3)node(fin1){};
		\draw[loosely dashed] (aux1.center)--(a-|aux1)node(x1)[anchor=north]{$y_1$};
		\def\ecartvert{1.5}	
		\draw (d.south west)++(0,-\ecartvert)node(d1){};
		\draw [->] (d1.center) -- coordinate [very near end] (f1)  +(0, 1) ; 
		\draw [->] (d1.center) -- +(5, 0); 
		\draw (d1)node(cons1)[anchor=south east,font=\fontsize{8}{144}\selectfont]{Consumption, Case 1};
		\draw [fill=red!40]  (f1) -- (f1 -| x1) coordinate (g1) -- (g1 |- d1) -- (d1.center) -- cycle ; 
		\draw [fill=red!40]  (f1-|b1) -- (f1 -| fin1) coordinate (g1) -- (g1 |- d1) -- (d1-|b1) -- cycle ; 
		\draw[<->](cons1 -| x1)--node[anchor=south]{$t_{\text{red}}$}(cons1-|b1);
		\draw[loosely dashed](x1)--(f1-|x1);
		\draw[loosely dashed](b1)--(f1-|b1);
		\draw[loosely dashed](fin1)--(f1-|fin1);
	}
	\color{blue}{
		\draw ([yshift=4cm] a)++(1.1,-1.65) node (aux2){};
		\draw [ultra thick,densely dashed] (aux2.center)--++(1.65,1.65)node(b2){}--++(2,-3)node(fin2){};
		\draw[loosely dashed] (aux2.center)--(a-|aux2)node(x2)[anchor=north]{$y_2$};
		\def\ecartvert{1.5}	
		\draw (d1)++(0,-\ecartvert)node(d2){};
		\draw [->] (d2.center) -- coordinate [very near end] (f2)  +(0, 1) ; 
		\draw [->] (d2.center) -- +(5, 0); 
		\draw (d2)node(cons2)[anchor=south east,font=\fontsize{8}{144}\selectfont]{Consumption, Case 2};
		\draw [fill=blue!40]  (f2) -- (f2 -| x2) coordinate (g2) -- (g2 |- d2) -- (d2.center) -- cycle ; 
		\draw [fill=blue!40]  (f2-|b2) -- (f2 -| fin2) coordinate (g2) -- (g2 |- d2) -- (d2-|b2) -- cycle ; 
		\draw[<->](cons2 -| x2)--node[anchor=south]{$t_{\text{red}}$}(cons2-|motorOff);
		\draw[loosely dashed](x2)--(f2-|x2);
		\draw[loosely dashed](b2)--(f2-|b2);
		\draw[loosely dashed](fin2)--(f2-|fin2);
		\draw[loosely dashed](d-|motorOff)--(d2-|motorOff);
	}
	\end{tikzpicture}
	}
	\caption{Two examples of consumption reduction for a cooling appliance. The solid line represents the evolution of the temperature without any reduction request; and the dashed lines the temperature evolution if a request comes when the appliance is at position $y_1$ or $y_2$ of its cycle. In the latter case it stops consuming upon receiving, until the temperature hits $T_{\max}$. Instantaneous consumptions and the reduction duration $t_{\text{red}}$ are displayed at the bottom.} 
	\label{fig:effacement_vs_cycle} 
\end{figure}
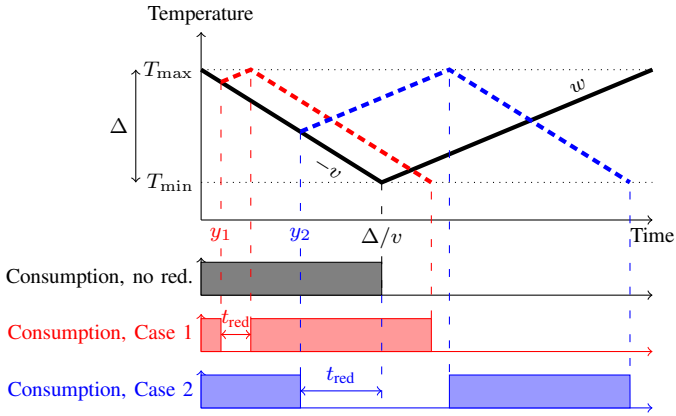

The conditions above are illustrated in Figure~\ref{fig:effacement_vs_cycle}, where we display the evolution of the temperature with time in several cases: \emph{(i)} if there is no request (solid line) and \emph{(ii)} if the appliance stops its temperature-modifier when at position $y$ in its cycle (two different values of $y$ are shown). 

Let us consider an appliance, responding to a reduction request by switching to OFF state. The consumption reduction is null if $y\geq \Delta/v$ (the appliance is already in OFF state), and otherwise it ends when one of the following events occur:
\begin{itemize}
	\item The temperature reaches the limit temperature and the system has to be turned ON again (Case $1$ in Figure~\ref{fig:effacement_vs_cycle}), i.e., after $\delta_T/w=\frac{yv}{w}$  ;	
	\item The normal cycle (without consumption reduction) would have ended and the cooling system would have turned off (Case $2$ in Figure~\ref{fig:effacement_vs_cycle}), which occurs after a duration $\Delta/v-y$.
\end{itemize}
Summarizing, the reduction duration $t_{\text{red}}$ is then
\begin{equation}\label{eq:duree_eff}
t_{\text{red}}=\min\left(yv/w, [\Delta/v-y]^+\right),
\end{equation}

\short{Consider an appliance}{Since appliances with $y\geq \Delta/v$ cannot provide any reduction (they are in OFF state), we focus on appliances} with $y<\Delta/v$.
\short{From~\eqref{eq:duree_eff},}{The first term being increasing in $y$ and the second one decreasing, the reduction duration is maximal \short{}{when both are equal, i.e., }when 
$
y= \frac{\Delta}{v} \times \frac{w}{v+w},
$
hence} the maximum reduction duration\short{ equals }{} 
$t_{\text{red}}^{\max}\define\frac{\Delta}{v+w}$.

\begin{figure}[htbp] 
	\centering 
{\footnotesize
	\begin{tikzpicture}[xscale=1.4,yscale=0.8]
	\coordinate (a) at (0,0);
	\draw [->] (a) -- +(0, 3.3)  node [above] {$t_{\text{red}}$}; 
	\draw [->] (a) -- +(5.5, 0) node [right] {$y$};	
	\draw [dashed] ([xshift=1.5cm] a) coordinate (seuil) -- +(0, 3) coordinate (effmax);
	\draw (seuil)node[anchor=north]{$\frac{\Delta}{v}\frac{w}{v+w}$};
	\draw[very thick,black!75] (a) -- node [above, rotate=48,near start] {slope $v/w$}(effmax) -- node [above, rotate=-29,near end] {slope $-1$} ++(3, -3) node [below] {$\Delta/v$}--++(0.8,0); 
	\draw[dashed] (effmax)--(a|-effmax)node[anchor=east]{$\frac{\Delta}{v+w}$};
	\draw [very thick,densely dashed, latex-latex] (0.8,1.64)coordinate(gauche) --node[anchor=north,fill=white,opacity=.8,text opacity=1]{$\left(\!1\!+\!\frac{w}{v}\!\right)\!\!\left[\!\frac{\Delta}{v+w}\!-\!t\right]^+$} +(2.1, 0)coordinate(droite); 
	\draw[dashed](gauche)--(a|-gauche)node[anchor=east]{$t$};
	\draw[dashed](gauche)--(a-|gauche)node[below]{$t\frac{w}{v}$};
	\draw[dashed](droite)--(a-|droite)node[below]{$\frac{\Delta}{v}-t$};
	\end{tikzpicture}
}
	\caption{Duration of the possible consumption reduction for an appliance \emph{versus} appliance cycle position $y$ upon receiving the request. The horizontal arrow identifies appliances able to provide a reduction during at least $t$.} 
	\label{fig:duree_effacement_vs_x} 
\end{figure}
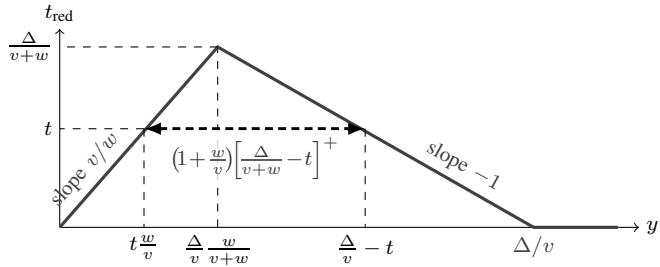

The reduction duration that a position-$y$ appliance can individually offer is plotted in Figure~\ref{fig:duree_effacement_vs_x}, with the range of appliances able to reduce during at least $t$.
\short{I}{From~\eqref{eq:duree_eff}, i}n the steady-state situation
the probability that an appliance can provide a reduction duration above $t$ is (see Figure~\ref{fig:duree_effacement_vs_x})
\begin{align}
{\mathbb P}(t_{\text{red}}>t)&= \frac{1}{\Delta(1/v+1/w)}\times \left(1+\frac wv\right)\left[\frac{\Delta}{v+w}-t\right]^+ \nonumber\\ 
&= \frac{w}{v+w}\left[1-t\frac{v+w}{\Delta}\right]^+.
\label{eq:distrib_duree_eff}
\end{align}

This corresponds to an average power reduction per appliance of $\frac{Pw}{v+w}\left[1-t\frac{v+w}{\Delta}\right]^+$, which when compared to the steady-state consumption $\frac{Pw}{v+w}$ gives us the following result.
\begin{proposition}\label{prop:schOne}
Over a duration $t$, the \schOne{} mechanism allows a relative power reduction of
\begin{equation}\label{eq:RedSchOne}
R_{\textrm{\schOne}} = \left[1-t\frac{v+w}{\Delta}\right]^+.
\end{equation}
\end{proposition}

\short{}{
That value is plotted in Figure~\ref{fig:power}, and compared to the upper bound of Proposition~\ref{prop:upper_bound}. As expected, there seems to be some room for offering larger reductions (although the upper bound does not consider the constant-power constraint). We investigate an improved version in Section~\ref{sec:schTwo}.
\begin{figure}[htbp]
{\footnotesize
	\begin{tikzpicture}[scale=3]
	
	\draw[->] (0,0) -- (2.3,0) node[right] {\emph{t}};
	\draw[->] (0,0) -- (0,1.3) node[above,font=\fontsize{8}{144}\selectfont] {Relative reduction};
	
	\draw[red,ultra thick,domain=0:1,smooth] plot (\x,{1-0.5*(\x)});
	\draw[red,ultra thick,domain=1:2.2,smooth] plot (\x,{0.5/(\x)});
	\draw[ultra thick,domain=0:1/1.2,smooth] plot (\x,{1-1.2*(\x)});
		
	\draw[dashed] (1,0.5)--(0,0.5);
	\draw[dashed] (1,0.5)--(1,0) ;
	\node[] at (1,-0.1) {$\frac{\Delta}{w}$};	
	\node[] at (-0.1,1) {$1$};
	\node[] at (-0.1,0.5) {$0.5$};
	\node[] at (-0.1,0) {$0$};
	\node[] at (0.8,-0.1) {$\frac{\Delta}{v+w}$};
	\node[] at (.38,.3) {$-\frac{v+w}{\Delta}$};
	\node[] at (.5,.84) {$-\frac{w}{2\Delta}$};
	
	\draw[red,ultra thick] (0.7,1.1)--(0.8,1.1) node[right,black,font=\fontsize{8}{144}\selectfont] {Upper bound}; 
	\draw[ultra thick] (0.7,0.95)--(0.8,0.95) node[right,black,font=\fontsize{8}{144}\selectfont] {\schOne};
	
	\fill [gray!50,nearly transparent, domain=1:2.2, variable=\x] (1,0)-- plot ({\x}, {0.5/\x})-- (2.2, 0.5/2.2)--(2.2,0)--cycle;
	
	\fill [gray!50,nearly transparent] (0,1)--(1,0.5)--(1,0)--(1/1.2,0)--cycle;
	
	\end{tikzpicture}
	}
	\caption{Reduced power (in proportion of the steady-state consumption $NP\frac{w}{v+w}$): upper bound and performance of \schOne{}. The shaded area illustrates the possible room for improvement.} 
	\label{fig:power} 
\end{figure}
}

\subsection{Implementing \schOne{} in practice}

Consider that the grid issues a request for a reduction of amplitude $A$ over a duration $t$. The aggregator can then directly use Proposition~\ref{prop:schOne} to know whether it can satisfy the request using \schOne{}: indeed it knows the number $N$ of appliances, thus over $t$ it can offer a reduction amplitude of $\frac{NPw}{v+w}R_{\textrm{\schOne}}$. 

Hence if $A\leq \frac{NPw}{v+w}\left[1-t\frac{v+w}{\Delta}\right]^+$ the aggregator can simply satisfy the request by broadcasting to all appliances a ``pulse'' request message interpreted as
\begin{quote}
\emph{``If the time $y$ since the beginning of your cycle is such that $\min\left(yv/w, \left[\Delta/y-x\right]^+\right)\geq t$, then turn off your engine as long as you can''},
\end{quote}
or equivalently 
\begin{quote}
\emph{``If you can reduce your demand immediately during at least $t$, do it''}.
\end{quote}

Note however that if the amplitude $A$ is strictly smaller than $N\frac{Pw}{v+w}R_{\textrm{\schOne}}$ then the reduction amplitude will exceed the one requested. If the aggregator wants to exactly offer a reduction amplitude $A$, we envision two simple possibilities:

\renewcommand{\theenumi}{\emph{\alph{enumi}}}
\begin{enumerate}
\item \label{it:approachOne}The reduction load can be taken by the appliances that can offer the longest reduction: in practice the aggregator would artificially increase the reduction duration so that the maximum reduction amplitude exactly equals $A$, so that the broadcasted message would be
\begin{quote}
\emph{``if you can reduce your demand immediately during at least $t'$, do it''},
\end{quote}
where the aggregator sets $t'\define\Delta\left(\frac{1}{v+w}-\frac{A}{NPw}\right)$.
\item \label{it:approachTwo}Alternatively, the reduction load can be shared evenly among all mobilizable appliances, i.e., the message broadcasted would be
\begin{quote}
\emph{``if you can reduce your demand immediately during at least $t$, do it with probability $p$, otherwise ignore this message''},
\end{quote}
where $p=\frac{A}{NP}\frac{v+w}{w R_{\textrm{\schOne}}}$ is computed by the aggregator as the ratio between the amplitude asked and the maximum possible amplitude.
\end{enumerate}
Of course other solutions are possible: we do not develop them here since we focus on providing the largest amplitudes.

\subsection{Aftermath of an \schOne{} reduction request}

We investigate here what happens just after a reduction request is satisfied.
For sake of clarity, we take the first approach (case~\ref{it:approachOne}) of the previous subsection to satisfy a request: all appliances that can offer a reduction of duration at least $t$ stop their temperature-modifying engine, whether that duration is the one asked or a strictly larger duration. Hence this is equivalent to offering a maximum reduction with \schOne{} for a duration $t$ (again, even if the actual duration asked is below $t$). Note however that the second approach (case~\ref{it:approachTwo}) could also be considered without major  difficulty\short{}{; mainly, the graphical interpretation we propose would be harder to interpret since it would involve probabilities (only a proportion of the highlighted appliances would be concerned by our comments)}.

In what follows, we study the operation of each appliance as a function of both $y$ (the appliance situation in its cycle upon receiving the reduction request) and $x$ (the time since the request was received).
As illustrated in Figure \ref{fig:duree_effacement_vs_x}, all appliances with $y$ in $[t\frac{w}{v}, \frac{\Delta}{v}-t]$ would be turned OFF. After this unique reaction to the request, appliances follow their usual functioning algorithm, i.e., remain OFF until the temperature hits the limit and then switch to ON, as shown in Figure~\ref{fig:effacement_vs_cycle}.

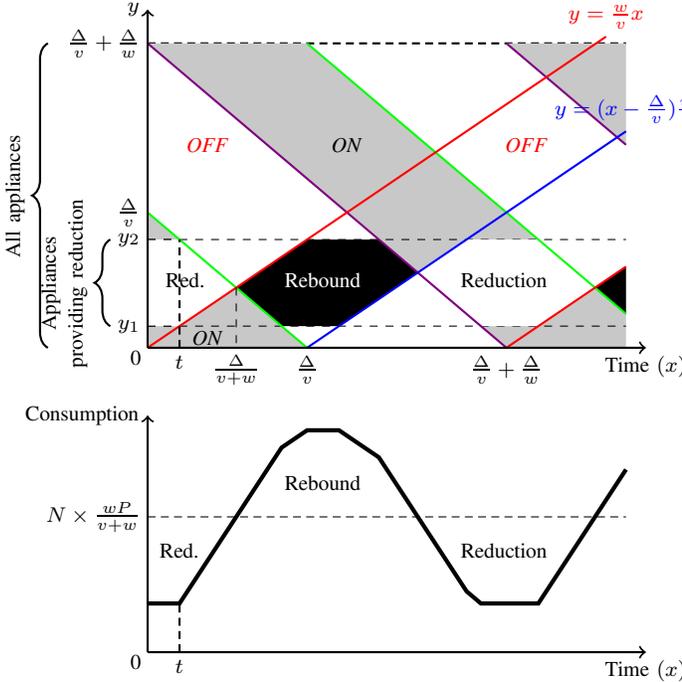
\begin{figure} [htbp]
\centering
{\footnotesize
\begin{tikzpicture}[xscale=.53,yscale=.45]

\gdef\stopx{12}

\draw[thick,densely dashed] (0,9)--(\stopx,9);

\fill[colorOn]  (0,9) -- (4,9)--(\stopx,13-\stopx)--(\stopx,0)--(9,0)--cycle;
\fill[colorOn]  (9,9) -- (\stopx,9)--(\stopx,18-\stopx)--cycle;
\fill[colorOn]  (0,0) -- (0,4)--(4,0)--cycle;

\draw[dashed] (0,3.2)--(\stopx,3.2);
\draw[dashed] (0,0.64)--(\stopx,0.64);

\node[red] at (1.5,6) {\emph{OFF}};\node[red] at (9.5,6) {\emph{OFF}};
\node at (1.5,0.3) {\emph{ON}};\node at (5,6){\emph{ON}};

\fill[white,draw opacity=0]  (0,3.2)--(0.8,3.2) -- (4/1.8,4/1.8*0.8) -- (0.8,0.64)--(0,0.64)-- cycle;
\fill[pattern=horizontal lines,draw opacity=0]  (0,3.2)--(0.8,3.2) -- (4/1.8,4/1.8*0.8) -- (0.8,0.64)--(0,0.64)-- cycle;
\node[] at (0.95,2) {Red.};
\fill[white,draw opacity=0]  (8,3.2)--(0.8+9,3.2) -- (4/1.8+9,4/1.8*0.8) -- (0.8+9,0.64)--(9-0.64,0.64)--(12.2/1.8,9-12.2/1.8)-- cycle;
\fill[pattern=horizontal lines,draw opacity=0]  (8,3.2)--(0.8+9,3.2) -- (4/1.8+9,4/1.8*0.8) -- (0.8+9,0.64)--(9-0.64,0.64)--(12.2/1.8,9-12.2/1.8)-- cycle;
\node[] at (8.95,2) {Reduction};

\fill[black]  (4/1.8,4/1.8*0.8) -- (4,3.2)--(5.8,3.2)--(12.2/1.8,9-12.2/1.8)--(4.8,0.64)--(3.36,0.64)--cycle; 
\node[white] at (4.4,2) {Rebound};
\fill[black]  (4/1.8+9,4/1.8*0.8) -- (\stopx,4/5*\stopx-4/5*9)--(\stopx,13-\stopx)--cycle; 

\draw[thick,green] (0,4) node[left,black] {$\frac{\Delta}{v}$}--(4,0) node [below,black] {$\frac{\Delta}{v}$};
\draw[thick,green] (4,9)--(\stopx,13-\stopx);
\draw[thick,violet] (0,9) node [left,black] {$\frac{\Delta}{v}+\frac{\Delta}{w}$} --(9,0) node [below,black]{$\frac{\Delta}{v}+\frac{\Delta}{w}$};
\draw[thick,violet] (9,9) --(\stopx,18-\stopx);
\draw[thick,domain=0:11.5,red] plot (\x,{0.8*(\x)}) node [above] {$y=\frac{w}{v}x$};
\draw[thick,domain=9:\stopx,red] plot (\x,{0.8*(\x-9)});
\draw[thick,domain=4:\stopx,blue] plot (\x,{0.8*(\x-4)}) node [above] {$y=(x-\frac{\Delta}{v})\frac{w}{v}$};

\draw[dashed] (4/1.8,4/1.8*0.8)--(4/1.8,0) node [below] {$\frac{\Delta}{v+w}$};
\draw[thick,densely dashed] (0.8,0) node[below] {$t$}--(0.8,3.2);
\draw[dashed] (0,3.2) node(y2)[left]{$y_2$}--(0.8,3.2);
\draw[dashed] (0,0.64) node(y1)[left]{$y_1$}--(0.8,0.64);
\node[] at (-0.3,-0.3) {0};

\draw[thick,->] (0,0) -- (12.5,0) node[below](rightend) {Time $(x)$};
\draw[thick,->] (0,0) -- (0,10) node[left] {$y$};

\draw [thick,decorate,decoration={brace,amplitude=6pt}](y1.west) -- (y2.west) node(provreq) [xshift=-.2cm,black,midway,anchor=south,rotate=90] {providing reduction};
\draw (provreq.north)++(.2,0)node[anchor=south,rotate=90]{Appliances};
\draw [thick,decorate,decoration={brace,amplitude=6pt}](-2.5,0) -- (-2.5,9) node [xshift=-.2cm,black,midway,anchor=south,rotate=90] {All appliances};


\gdef\ordonnee{-9}
\gdef\hauteur{7}
\FPeval{ecart}{3.2-0.64}
\gdef\const{4}	

	\draw[thick,->] (0,\ordonnee) -- (0,\ordonnee-|rightend) node[below](rightendbas) {Time $(x)$};
	\draw[thick,->] (0,\ordonnee) -- (0,\hauteur+\ordonnee) node[left] {Consumption};	
	\node[] at (-0.3,-0.3+\ordonnee) {0};
	
	\draw[densely dashed] (0,\ordonnee+\const)node[left]{$N\times\frac{wP}{v+w}$}--(\stopx,\ordonnee+\const);
	
	\FPeval{reducmax}{\const-\ecart+\ordonnee}
	\FPeval{consoAux}{\const+\ecart-(4-3.36)*4/5+\ordonnee}
	\draw[ultra thick](0,\reducmax)--(0.8,\reducmax)--(4/1.8,\ordonnee+\const)--(3.36,\consoAux)--(4,\ordonnee+\const+\ecart)--(4.8,\ordonnee+\const+\ecart)--(5.8,\ordonnee+\const+\ecart-4/5)--(12.2/1.8,\ordonnee+\const)--(8,\reducmax+0.36)--(9-0.64,\reducmax)--(0.8+9,\reducmax)--(4/1.8+9,\ordonnee+\const)--++(\stopx-4/1.8-9,9/5*\stopx-9/5*4/1.8-9/5*9)
	;
	\draw[thick,densely dashed] (0.8,\ordonnee) node[below] {$t$}--(0.8,\reducmax);

	\node at (4.4,\const+\ordonnee+1) {Rebound};
	\node at (8.95,\const+\ordonnee-1) {Reduction};
	\node at (0.8,\const+\ordonnee-1) {Red.};
	
	\end{tikzpicture}
	}
	\caption{\emph{(Top)} Appliance states \emph{vs} time, after an \schOne{} duration-$t$ reduction request. The appliances affected are those in the $[y_1,y_2]$ range. Horizontal line patterns indicate appliances that are OFF while they would be ON without the request (hence a reduction), black zones indicate the opposite (hence a rebound: larger consumption than without the request). The other zones are ``neutral'': the appliance is in the state it would be without the request.\newline
	\emph{(Bottom)} The resulting consumption pattern just after the request, reflecting the height of reduction and rebound zones. The analytical expression (piecewise linear) can easily be obtained from the linear functions of the top part. The dashed line is the consumption without the request.} 
	\label{fig:aftermath_simple} 
\end{figure}

This behavior is the simplest possible, and leads to a modified ON-OFF pattern and a new consumption curve, illustrated in Figure~\ref{fig:aftermath_simple}.
\short{}{The figure displays the appliances affected by the request, and the consequences on the global consumption. }By design, we have a constant consumption during $t$, corresponding to the requested reduction. Then consumption increases linearly and reaches the steady-state consumption some time $\frac{\Delta}{v+w}$ after the request; we then enter a \emph{rebound} phase where consumption exceeds the steady-state one.

In the next section, we exploit that behavior to extend the scheme so as to offer larger--or longer--reductions\short{}{, thereby getting closer to the upper bound plotted in Figure~\ref{fig:power}}.

\section{\schTwo{}: Extending the reduction by coordinating appliances} \label{sec:schTwo}

The \schOne{} mechanism relies on appliances individually performing a consumption reduction during the asked duration $t$. In this section, we suggest to coordinate appliances so that some reduce their consumption \emph{after} the start of the reduction, in order to \emph{compensate} for the limited reduction durations of others. Like for \schOne{}, the whole reduction can be triggered by a single broadcast message, but not all appliances react at the same time to that message.

\short{}{Our reasoning will use the graph of individual appliance states over time of Figure~\ref{fig:aftermath_simple}.}

\subsection{Principle} 

The idea is to start like in \schOne{}, but to additionally have new appliances contribute when those initially providing the reduction stop reducing. In Figure~\ref{fig:aftermath_simple} this occurs after $t$, but \schTwo{} will allow to provide longer reductions. Hence we will denote by $\tilde t$ the time when the first appliances stop reducing, and by $t$ the total reduction duration.

\begin{figure} [htbp]
{\footnotesize
	\centering
	\begin{tikzpicture}[yscale=.6]

	\fill[colorOn]  (0,0) -- (0,5)--(5,0)--cycle;
	\fill[colorOn]  (0,9) -- (5,9)--(6,8)--(6,3)--cycle;

	\node (A1)  at (0,5-1.5) {\hpat{$A_1$}};
	\node (A2)  at (1.5,5-1.5) {\hpat{$A_2$}}; 
	\node (A3)  at (5/1.8,5/1.8*0.8) {\hpat{$A_3$}};
	\node (A4)  at (1.5,1.5*0.8) {\hpat{$A_4$}};
	\node (A5)  at (0,1.5*0.8) {\hpat{$A_5$}};
	\fill[white,draw opacity=0]  (A1.center)node[black,left]{$y_2$}--(A2.center) -- (A3.center) -- (A4.center)--(A5.center)node[black,left]{$y_1$}-- cycle;
	\fill[pattern=horizontal lines,draw opacity=0]  (A1.center)--(A2.center) -- (A3.center) -- (A4.center)--(A5.center)-- cycle;
	\node[fill=white,opacity=.5,text opacity=1] at (0.95,2.5) {Red.};

	\node (B1)  at (A4) {\hpat{B1}};
	\node (B3)  at (3.6,0) {\hpat{B3}};
	\node (B2)  at (B1-|B3) {\hpat{B2}};
	\node (B4)  at (1.5*0.8/1.8+1.5,0){\hpat{B4}};
	\node (B5p)  at (3.6,-2.5){\hpat{B5'}};
	\node (B3p)  at (3.6,9) {\hpat{B3'}};
	\node (B4p)  at (1.5*0.8/1.8+1.5,9) {\hpat{B4'}};
	\node (B5)  at (3.6,6.42) {\hpat{B5}};
	\fill[white,draw opacity=0]  (B1.center)--(B2.center) -- (B3.center) -- (B4.center)-- cycle;
	\fill[pattern=vertical lines,draw opacity=0]  (B1.center)--(B2.center) -- (B3.center) -- (B4.center)-- cycle;
	\node[fill=white,opacity=.5,text opacity=1] at (2.8,.5) {Red.};
	\fill[white,draw opacity=0]  (B3p.center)--(B4p.center) -- (B5.center)-- cycle;
	\fill[pattern=vertical lines,draw opacity=0]  (B3p.center)--(B4p.center) -- (B5.center)-- cycle;
	\node[fill=white,opacity=.5,text opacity=1] at (3,8.4) {Red.};
	\draw[dotted] (B5.center)--(B5-|0,0) node[left]{$y_3$};

	\fill[black] (A3.center)--(B5|-0,5/1.8*0.8-3.6+5/1.8)--(B5|-0,5/1.8*0.8+4/5*3.6-4/1.8)--cycle;
	
	\draw[thick,green] (0,5) node[left,black] {$\frac{\Delta}{v}$}--(5,0) node [below,black] {$\frac{\Delta}{v}$};
	\draw[thick,green] (5,9) --(6,8) node [below,black] {};
	\draw[thick,violet] (0,9) node [left,black]{$\frac{\Delta}{v}\!\!+\!\!\frac{\Delta}{w}\!$}--(6,3);
	\draw[thick,dashed,violet] (0,0) --(3.7,-3.7) ;

	\draw [thick,decorate,decoration={brace,amplitude=6pt}](-1.4,0) -- (-1.4,9) node [xshift=-.2cm,black,midway,anchor=south,rotate=90] {All appliances};
	\draw [thick,decorate,decoration={brace,amplitude=6pt}] (-.6,1.5*0.8)--(-.6,5-1.5) node [xshift=-.2cm,black,midway,anchor=south,rotate=90] {First batch};
	\draw [thick,decorate,decoration={brace,amplitude=6pt}] (-.8,0|-B5p)--(-.8,1.5*0.8-.05) node [xshift=-.2cm,black,midway,anchor=south,rotate=90] {Second batch};
	\draw [white, line width=1mm, densely dotted]  (-.9,0|-B5p)--(-.9,0);
	\draw [thick,decorate,decoration={brace,amplitude=6pt}] (-.8,0|-B5)--(-.8,10.2) node [xshift=-.2cm,black,midway,anchor=south,rotate=90] {Second batch};
	\draw [white, line width=1mm, densely dotted]  (-.9,10.2)--(-.9,9);

	\draw[thick,domain=0:6,smooth,red] plot (\x,{0.8*(\x)}) node [above,black] {$y=\frac{w}{v}x$};
	\draw[thick,domain=5:6,smooth,blue] plot (\x,{0.8*(\x-5)}) node (auxdroite) [above,black] {$y=\frac{w}{v}(x-\frac{\Delta}{v})$};
	\draw (0,9.2) node (auxhaut){};
	\draw (B5p) node [right,black] {$y=y_1-(1+\frac{w}{v})(t-\tilde{t})$};

	
	\node  at (1.4,-0.22) {$\tilde{t}$};
	
	\node at (3.8,-0.2) {$t$};
	
	\node[] at (1.5,-2)[rotate=-35] {slope $-1$};
	
	\draw[dashed] (0,0-|A2.center)--(1.5,10.2)--(3.6,10.2)--(B5p.center)--(A4.center)--(B2.center)
	(1.5,10.2)--(B5.center);

	
	\node[] at (1,0.3) {\emph{ON}};\node[] at (1,8.7) {\emph{ON}};
	\node[red] at (2.4,5) {\emph{OFF}};

	\draw (3.61,0) node (t){};
	\fill[white,opacity=.7] (t.center)--(t.center|-auxhaut.north)--(auxhaut.north-|auxdroite.east)--(auxdroite.east|-t.center)--cycle;
	
	\def\stopx{6.5}
	\draw[thick,dashed] (0,9)--(\stopx,9);
	\draw[thick,->] (0,0) -- (\stopx,0) node[below](time) {Time};
	\draw (time)++(0,-.1) node [below]{($x$)};
	\draw[thick,->] (0,0) -- (0,10.5) node[above] {$y$};	
	\node[] at (-0.2,0) {$0$};

\end{tikzpicture}
\caption{Appliance states over time for \schTwo{} up to the requested reduction duration $t$, to offer a constant reduction over $t$. Line-patterned areas represent consumption reductions, and the black area is a rebound (appliances are ON while they would be OFF without the reduction request).} 
\label{fig:schTwo}  
}
\end{figure}
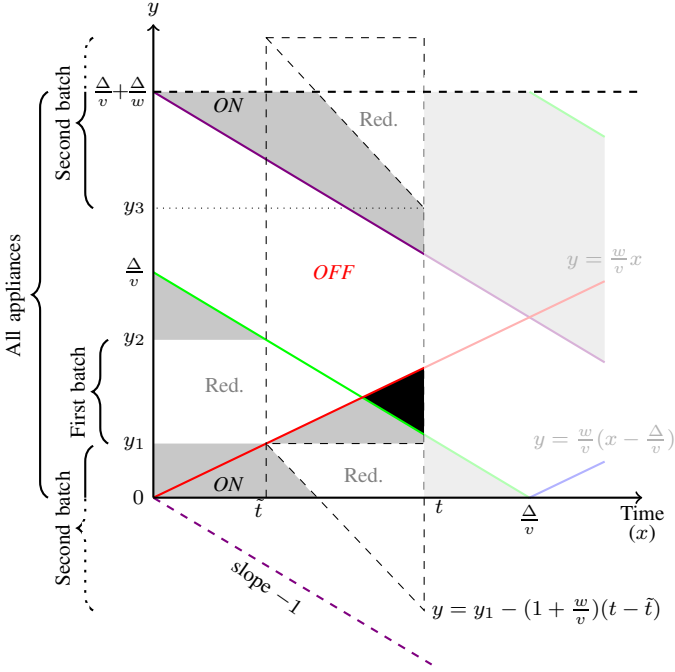

The functioning of \schTwo{} is depicted in Figure~\ref{fig:schTwo}: we start with a reduction as with \schOne{}, relying on a first batch of appliances. After time $\tilde t$, some appliances begin to stop reducing, either because they have to turn ON, or because they would have turned OFF without the request. Cumulating those two causes, the overall consumption increases at a constant speed $NP\left(1+\frac{w}{v}\right)$ until all the concerned appliances for either cause are affected, i.e., until time 
$
\short{}{\hat t\define }\min\left(\frac{\Delta}{v}-y_1,y_2\frac{v}{w}\right).
$
The addition of \schTwo{} with respect to \schOne{} is to involve a \emph{second batch of appliances}, entering the reduction \emph{gradually} from time $\tilde t$, exactly at the speed $N\left(1+\frac{w}{v}\right)$: doing so, the extra reduction compensates the consumption increase\short{}{ that starts at $\tilde t$}. \short{}{The goal is that each of those new appliances (individually) offer a constant reduction until time $t\leq \hat t$.}

This cannot be done infinitely, since we face two constraints:
\begin{enumerate}
\item Those newly involved appliances need to be able to offer a reduction until time $t$, thus
\begin{itemize}
\item they should be able to stay OFF from the moment they are supposed to participate until time $t$, and
\item without the request they would have been ON during that period.
\end{itemize}
\item Those appliances cannot be among those already involved in the first batch.
\end{enumerate}

The first batch involves appliances with cycle positions in the interval $[y_1,y_2]$ as shown in Figure \ref{fig:schTwo}, with a constant reduction until time $\tilde t$.
From that instant, appliances from the second batch enter progressively, at ``speed'' $1+w/v$: more precisely, between time $\tilde t$ and $\tilde t+x$ we need an extra proportion $(1+w/v)\frac{x}{\Delta/v+\Delta/w}=\frac{w}{\Delta}x$ of all appliances to contribute to the reduction. For this, we rely on appliances whose cycle position (upon emission of the request, and modulo $\Delta/v+\Delta/w$) is in the range $[y_1-(1+w/v)x,y_1]$, as illustrated in Figure~\ref{fig:schTwo}.

\subsection{Reduction potential with \schTwo}

The following proposition quantifies how much reduction can be offered for a given reduction duration $t$.

\begin{proposition}\label{prop:schTwo}
Over a duration $t\leq t_{\text{schTwo}}^{\max}\define \Delta\frac{v+2w}{(v+w)^2}$, the \schTwo{} mechanism allows a relative power reduction of
\begin{equation}\label{eq:RedSchTwo}
R_{\textrm{\schTwo}} = 1-t\frac{v+w}{v+2w}\frac{w}{\Delta}.
\end{equation}
\end{proposition}

\begin{IEEEproof}
The reduction duration $t$ is imposed, but not the duration $\tilde t$ of the reduction from the first batch, which also determines the values of $y_1$ and $y_2$ and thus the relative reduction amplitude $R$:
\begin{eqnarray}
R &=&\frac{y_2-y_1}{\Delta/v}
\ =\  \short{}{\frac{\Delta/v - (1+w/v)\tilde t}{\Delta/v}\nonumber\\
&=&} 1 - \frac{(v+w)\tilde t}{\Delta},\label{eq:amplitude_ttilde}
\end{eqnarray}
hence to maximize the reduction amplitude, one has to select $\tilde t$ as small as possible\short{}{, as Figure~\ref{fig:schTwo} also illustrates}.

On the other hand, the constraints highlighted above correspond mathematically, with the notations of Figure~\ref{fig:schTwo}, to
\begin{eqnarray}
y_3&\geq& y_2\label{eq:constr1}\\
y_3&\geq& \Delta/v+\Delta/w-t\label{eq:constr2}\\
y_1&\leq& \Delta/v-t\label{eq:constr3}\\
t-\tilde t&\leq& v/w(\tilde t+y_1).\label{eq:constr4}
\end{eqnarray}
Condition~\eqref{eq:constr1} states that the second batch is separated from the first batch; \eqref{eq:constr2} and~\eqref{eq:constr3} mean that the hashed triangle in Figure~\ref{fig:schTwo} is included in an ON zone, i.e., without the request all second-batch appliances would have been ON from their entering the service until $t$. Finally, \eqref{eq:constr4} indicates that when entering the service, any second-batch appliance is able to remain OFF until $t$ (i.e., its temperature stays within $[T_{\min},T_{\max}]$): the appliance entering the service at time $\tilde t+x$ has stayed ON during $\tilde t+y_1 - xw/v$, and can therefore remain OFF for $v/w$ times that duration. We need the result to exceed the remaining reduction time $t-(\tilde t+x)$, which yields
\[
t-\tilde t-x\leq v/w\left(\tilde t+y_1 - xw/v\right)\quad \forall x\in[0,t-\tilde t],
\]
and simplifies to~\eqref{eq:constr4}.

Now, we use the expressions of $y_1$, $y_2$, and $y_3$
\begin{equation}\label{eq:y_schTwo}
\left\{
\begin{array}{rcl}
y_1 &=&\tilde t w/v\\
y_2 &=&\Delta/v - \tilde t\\
y_3 &=&\Delta/v + \Delta/w + \tilde t w/v - (1+w/v)(t-\tilde t)
\end{array}
\right.
\end{equation}
to rewrite~\eqref{eq:constr1}-\short{\eqref{eq:constr3}}{\eqref{eq:constr4}} in terms of the decision variable $\tilde t$, and the known values $\Delta,v,w,$ and $t$. \short{We respectively obtain}{Conditions~\eqref{eq:constr1},~\eqref{eq:constr2}, and~\eqref{eq:constr3} respectively give}
\begin{eqnarray}
\tilde t&\geq&\frac12\left( t-\frac{\Delta}{w+w^2/v}\right)\label{eq:constr1prime}
\\
\tilde t&\geq& \frac{w}{v+2w}t\label{eq:constr2prime}
\\
\tilde t&\leq& (\Delta-vt)/w,\label{eq:constr3prime}
\end{eqnarray}
while~\eqref{eq:constr4} exactly gives~\eqref{eq:constr2prime} and is therefore redundant.

But~\eqref{eq:constr1prime} is also redundant. Indeed,~\eqref{eq:constr2prime} and~\eqref{eq:constr3prime} give $\frac{w}{v+2w}t\leq\tilde t\leq \frac{\Delta}{w}-\frac{v}{w}t$, yielding $\frac{t}{v+2w}\leq\frac{\Delta}{(v+w)^2}$. And~\eqref{eq:constr2prime} can be rewritten as
$\tilde t\geq \frac12\left(t-\frac{v}{v+2w}t\right)$; then plugging the previous inequality gives
$\tilde t\geq \frac12\left(t-\frac{v\Delta}{(v+w)^2}\right)
= \frac12\left(t-\frac{\Delta}{v+2w+w^2/v}\right)$,
a condition stricter than~\eqref{eq:constr1prime}.

Summarizing, a duration-$t$ reduction is possible with \schTwo{} if and only if~\eqref{eq:constr2prime} and~\eqref{eq:constr3prime} are jointly satisfiable, i.e., if
$
t\leq \Delta\frac{v+2w}{(v+w)^2}$.
Under that condition, to maximize the reduction one must choose the smallest $\tilde t$, which is given in~\eqref{eq:constr2prime} as 
\begin{equation}\label{eq:ttilde_schTwo}
\tilde t = \frac{w}{v+2w}t,
\end{equation}
yielding a relative reduction, from~\eqref{eq:amplitude_ttilde}, of
$
R \!=\! 1 - \frac{v+w}{v+2w}\frac{wt}{\Delta}.\!
$
\short{}{This establishes the proposition.}
\end{IEEEproof}
Comparing with Proposition~\ref{prop:schOne}, we remark that \schTwo{} allows longer reductions than \schOne. However, note that~\eqref{eq:RedSchTwo} is strictly positive for $t=t_{\text{\schTwo}}^{\max}$, suggesting that some reduction is possible for a duration larger than $t_{\text{\schTwo}}^{\max}$. However this is not doable in the strict (and quite simple) sense with which we defined \schTwo: another combination of the appliances is needed.

\subsection{Implementing \schTwo{} from a single broadcast message}

As for \schOne, we can implement \schTwo{} by broadcasting a single message to all appliances. Here, to obtain a maximum amplitude the manager could \short{simply send the request duration $t$: each appliance would then compute $y_1$ and $y_2$  from~\eqref{eq:y_schTwo} and~\eqref{eq:ttilde_schTwo}:
\[
y_1 = \frac{w^2}{v(v+2w)}t\qquad ;\qquad
y_2 = \frac{\Delta}{v}-\frac{w}{v+2w}t.
\] and interpret the message as}{compute the values of $y_1,y_2$ shown in Figure~\ref{fig:schTwo}, and send them together with $t$ so that each appliance interprets the message as}
\begin{quote}
\emph{``If the time $y$ since the beginning of your cycle is such that $y\in[y_1,y_2]$, then turn off your engine as long as you can.\\
Otherwise, if there is an $x\in[0,t-y_1\frac{v}{w}]$ such that \[y \equiv y_1-(1+\frac{w}{v})x \textrm{ mod } (\Delta/v+\Delta/w),\] then wait until time $y_1\frac{v}{w}+x$ to turn off your engine as long as you can.''}
\end{quote}
\short{}{Note that assuming the appliances know $v$ and $w$, sending $\tilde t$ is unnecessary since it can be computed as $\tilde t=y_1\frac{v}{w}$.

Also, the manager could simply send the request duration $t$, and each appliance would behave as explained before, computing $y_1$ and $y_2$ themselves from~\eqref{eq:y_schTwo} and~\eqref{eq:ttilde_schTwo}:
\[
y_1 = \frac{w^2}{v(v+2w)}t\qquad ;\qquad
y_2 = \frac{\Delta}{v}-\frac{w}{v+2w}t.
\]
}

This method would provide the maximum \schTwo{} reduction amplitude possible, which we can denote by $A$. But as with \schOne{}, smaller amplitudes $A'<A$ can also be offered simply by having each appliance obey the message with probability $A'/A$ and ignore it otherwise. That probability would then be added to the broadcasted message.

\section{Discussion}

We discuss here the applicability of our schemes, their performance, and some variants and directions for future work.

\subsection{Applicability of the reduction schemes}

The mechanisms \schOne{} and \schTwo{} are both very simple, involving a simple calculation and at most one action from each appliance (turn OFF at a specific instant).\short{}{ Hence no computational power is needed on the appliance side.}

Moreover, in terms of communications our mechanisms are extremely lightweight: a reduction request (which should occur quite rarely) only involves the broadcast of one single message, containing very little information\short{: the reduction duration, plus possibly a ``probability to participate''.}{. As we saw, sending the reduction duration is sufficient to get the maximum amplitude, and one can possibly add a ``probability to participate'' field to provide smaller amplitudes.}

Hence we think both mechanisms are quite easily implementable in the context of the Internet of Things, even with very limited computational and communication capabilities. \short{}{Nevertheless, we point out here that the duty cycle constraints for connected objects~\cite{sheng2013survey} may incur some delays, which are ignored in our model but deserve further study. }

\subsection{Possible reductions with \schOne{} and \schTwo}

We compare here the performance results of Propositions~\ref{prop:upper_bound}, \ref{prop:schOne}, and~\ref{prop:schTwo}, in terms of the maximum reduction that can be offered over some duration $t$. That reduction (in proportion of the average consumption) is plotted in Figure~\ref{fig:power_new}\short{}{ for two different settings}. 
\gdef\CDelta{1}
\gdef\Cv{0.4}\gdef\Cvold{0.4}  
\gdef\Cw{1}
\gdef\xsup{2.2}
\begin{figure}[htbp]
{\footnotesize
\begin{tikzpicture}[scale=3]
\draw[->] (0,0) -- (2.3,0) node[below] {\begin{tabular}{c}Reduction\\duration \emph{t}\end{tabular}};
\draw[->] (0,0) -- (0,1.1) node[above] {\begin{tabular}{c}Relative\\reduction\end{tabular}};
\draw (0,1) node[left]{$1$};
\draw (0,0) node[left]{$0$};

\draw[red,ultra thick,domain=0:\CDelta/\Cw,smooth] plot (\x,{1-\x*\Cw/(2*\CDelta)});
\draw[red,ultra thick,domain=\CDelta/\Cw:\xsup,smooth] plot (\x,{\CDelta/(2*\Cw*\x)});
\node (aux_upper) at (\CDelta/\Cw,0.5) {};
\draw[red,dashed] (aux_upper-|0,0)node[left]{$0.5$}--(aux_upper.center)--(aux_upper|-0,0)node[below]{$\frac{\Delta}{w}$};

\FPeval{tmaxOne}{\CDelta/(\Cv+\Cw)}
\draw[ultra thick] (0,1)--(\tmaxOne,0) node[below]{$\frac{\Delta}{v+w}$}--(\xsup,0);

\FPeval{tmaxTwo}{\CDelta*(\Cv+2*\Cw)/(\Cv+\Cw)^2}
\FPeval{auxschTwo}{1-(\Cv+\Cw)/(\Cv+2*\Cw)*\Cw/\CDelta*\tmaxTwo}
\FPeval{tmaxTwoFalse}{\CDelta/\Cw*(\Cv+2*\Cw)/(\Cv+\Cw)}
\draw[blue,dashed] (\tmaxTwo,\auxschTwo) -- (\tmaxTwo,0) --++(0,-.15)  node[below] {$\frac{\Delta(v+2w)}{(v+w)^2}$};
\draw[blue,ultra thick] (0,1)--(\tmaxTwo,\auxschTwo);
\draw[blue,ultra thick] (\tmaxTwo,0)--(\xsup,0);
\draw[blue,very thick,densely dotted](0,1)--(\tmaxTwoFalse,0)node[below]{$\frac{\Delta(v+2w)}{w(v+w)}$};

\draw[red,ultra thick] (aux_upper|-0,1.3)coordinate(ref1)--++(0.1,0) node[right,black] {Upper bound}; 
\draw[ultra thick] (ref1)++(0,-.1)--++(0.1,0) node[right,black] {\schOne};
\draw[blue,ultra thick] (ref1)++(0,-.2)--++(0.1,0) node[right,black] {\schTwo};

\fill [gray!50,nearly transparent] (0,1)-- plot[domain=\CDelta/\Cw:\xsup, variable=\x] ({\x}, {\CDelta/(2*\Cw*\x)})-- (\xsup,0)--(\tmaxOne,0)--cycle;

\end{tikzpicture}
\short{}{
\gdef\Cv{2}
\begin{tikzpicture}[scale=3]
\draw[->] (0,0) -- (2.3,0) node[below] {\begin{tabular}{c}Reduction\\duration \emph{t}\end{tabular}};
\draw[->] (0,0) -- (0,1.3) node[above] {\begin{tabular}{c}Relative\\reduction\end{tabular}};
\draw (0,1) node[left]{$1$};
\draw (0,0) node[left]{$0$};

\draw[red,ultra thick,domain=0:\CDelta/\Cw,smooth] plot (\x,{1-\x*\Cw/(2*\CDelta)});
\draw[red,ultra thick,domain=\CDelta/\Cw:\xsup,smooth] plot (\x,{\CDelta/(2*\Cw*\x)});
\node (aux_upper) at (\CDelta/\Cw,0.5) {};
\draw[red,dashed] (aux_upper-|0,0)node[left]{$0.5$}--(aux_upper.center)--(aux_upper|-0,0)node[below]{$\frac{\Delta}{w}$};

\FPeval{tmaxOne}{\CDelta/(\Cv+\Cw)}
\draw[ultra thick] (0,1)--(\tmaxOne,0) node[below]{$\frac{\Delta}{v+w}$}--(\xsup,0);

\FPeval{tmaxTwo}{\CDelta*(\Cv+2*\Cw)/(\Cv+\Cw)^2}
\FPeval{auxschTwo}{1-(\Cv+\Cw)/(\Cv+2*\Cw)*\Cw/\CDelta*\tmaxTwo}
\FPeval{tmaxTwoFalse}{\CDelta/\Cw*(\Cv+2*\Cw)/(\Cv+\Cw)}
\draw[blue,dashed] (\tmaxTwo,\auxschTwo) -- (\tmaxTwo,0) --++(0,-.15) node[below] {$\frac{\Delta(v+2w)}{(v+w)^2}$};
\draw[blue,ultra thick] (0,1)--(\tmaxTwo,\auxschTwo);
\draw[blue,ultra thick] (\tmaxTwo,0)--(\xsup,0);
\draw[blue,very thick,densely dotted](0,1)--(\tmaxTwoFalse,0)node[below]{$\frac{\Delta(v+2w)}{w(v+w)}$};

\draw[red,ultra thick] (aux_upper|-0,1.2)coordinate(ref1)--++(0.1,0) node[right,black] {Upper bound}; 
\draw[ultra thick] (ref1)++(0,-.1)--++(0.1,0) node[right,black] {\schOne};
\draw[blue,ultra thick] (ref1)++(0,-.2)--++(0.1,0) node[right,black] {\schTwo};

\fill [gray!50,nearly transparent] (0,1)-- plot[domain=\CDelta/\Cw:\xsup, variable=\x] ({\x}, {\CDelta/(2*\Cw*\x)})-- (\xsup,0)--(\tmaxOne,0)--cycle;

\end{tikzpicture}
}
} 
\caption{Maximum relative reduction (in proportion of the aggregated average power consumption $NP\frac{w}{v+w}$) versus reduction duration, for $v=\Cvold$ \short{,}{\emph{(top)} and $v=\Cv$ \emph{(bottom)}, with} $\Delta=\CDelta$, $w=\Cw$.} 
\label{fig:power_new} 
\end{figure}
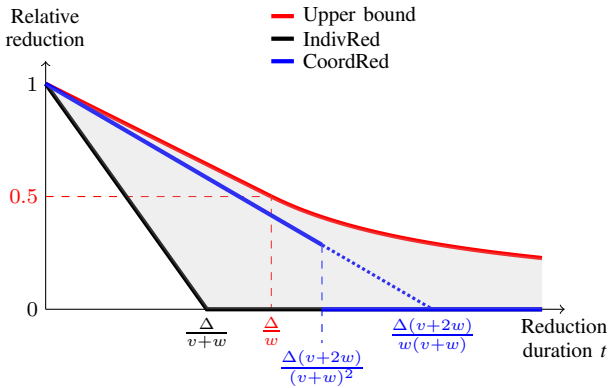

We observe that~\schTwo{} offers a considerable improvement with respect to \schOne{}, for $t\leq \frac{\Delta (v+2w)}{w(v+w)}$. Depending on the relative values of $v$ and $w$, this duration can be quite restrictive (it can be lower than the usual OFF duration $\Delta/w$ when $v$ is large). But in all cases, the maximum duration and the maximum amplitude are higher with~\schTwo.

\subsection{Mechanism variant to obtain demand increases}

This whole paper has been formulated in terms of demand \emph{reductions}, since the most frequent concern is about managing scarce energy production. But in a few occasions, and especially with renewable energies, we can have an over-production and want to temporarily \emph{increase} demand instead of decreasing it, as discussed in~\cite{gotz2013negative}. 

This is very easily doable within the context of this paper: just by exchanging the roles of $v$ and $w$ in Propositions~\ref{prop:upper_bound},~\ref{prop:schOne}, and~\ref{prop:schTwo}, we obtain the maximum increase in consumption, as a proportion of the \emph{average non-consumption}, that is $NP\frac{v}{v+w}$. 
More explicitly, adapting Equations~\eqref{eq:RedUpper_bound},~\eqref{eq:RedSchOne}, and~\eqref{eq:RedSchTwo} respectively yield that over a duration $t$:
\begin{enumerate}
\item One cannot get an average consumption increase of more than 
\[
\max\left(1 - \frac{vt}{2\Delta}, \frac{\Delta}{2vt}\right) \times NP\frac{v}{v+w}\text{ watts;}
\]
\item Adapting the \schOne{} mechanism allows a constant consumption increase of
\[
\left[1-t\frac{v+w}{\Delta}\right]^+ \times NP\frac{v}{v+w}\text{ watts;}
\]
\item If $t\leq\Delta\frac{w+2v}{(v+w)^2}$, adapting the \schTwo{} mechanism allows a constant consumption increase of 
\[
\left(1-t\frac{v+w}{2v+w}\frac{v}{\Delta}\right)\times NP\frac{v}{v+w}\text{ watts.}
\]
\end{enumerate}

\subsection{Possible extensions}

We discuss in this section some additional aspects that can be taken into account in future work\short{}{, to enrich our basic \schOne{} and \schTwo{} mechanisms}. 
\short{}{We evoke appliance heterogeneity, possible transmission issues for the reduction message, and more elaborate coordination schemes.}

\subsubsection{Managing several types of appliances}

Our model assumes all appliances are identical, with the same consumed power in ON state, the same temperature limits, and the same heating and cooling speeds. In practice, we will want to leverage the reduction potential of an heterogeneous set of appliances, with different parameters.

With our results formulated in terms of the reduction duration, it is quite simple to classify appliances into classes (appliances within a class being identical), so that the total reduction one can get over a time $t$ is just the sum of the reductions (in Watts) we can get from all classes.

One can also envision richer mechanisms, where classes are coordinated so that the reduction offered by each class is not of constant power, but the sum is. This \short{}{is beyond the scope of the paper, but }may be worth considering especially if heterogeneity among classes is large. 

\subsubsection{Coping with transmission errors and delays}

Our model ignores transmission issues, assuming that all appliances immediately receive the demand reduction broadcast message. 

In practice, problems such as losses and delays can occur\short{.}{, in particular within the Internet of Things (IoT) context.
}
Indeed, IoT protocols often involve some duty cycle constraints~\cite{sheng2013survey}, meaning that nodes cannot emit more than a given proportion of the time. Hence a node may have to wait before being allowed to forward a reduction request message. 
Also, those protocols~\cite{farrell2017LPWAN,lora2015lora,etsi2014LTN} are\short{}{ very lightweight and} subject to collisions, which incurs extra delays (due to retransmissions) or  message losses.

Those aspects should be considered when applying our mechanisms. The difficulties may be easily manageable (e.g., by sending the reduction request a bit ahead of time to absorb all possible delays, and by implementing reliability-oriented protocols), but they should not be forgotten.

\subsubsection{Combining more than two batches, allowing more complex appliance behavior}

Our \schOne{} and \schTwo{} schemes respectively involve one and two appliance batches to provide a reduction, and what we ask each appliance is extremely simple: ``\emph{switch to OFF state at this specific instant}''.

One can imagine more complex schemes, that would combine more batches and/or involve more subtle behaviors of individual appliances. This direction leaves some space for future works, especially to overcome the duration limitation of \schTwo. Nevertheless, this should make the analysis of those schemes more complex. Also, our schemes have the advantage of limiting the number of ON-OFF switches, each one possibly involving some energy costs (ignored in our model).

\section{Conclusion and perspectives} \label{sec:conclusion_perspectives}

This paper investigates how a large number of temperature-modifying appliances\short{}{ (e.g., water heaters, A/C or heating systems, fridges and freezers)}, when connected, offer new opportunities for demand flexibility. Based on a simple mathematical model, we quantify the power by which those appliances can reduce their aggregated consumption over a given period of time, while respecting individual temperature constraints.

In particular, we describe and analyze two mechanisms to coordinate appliances and offer significant power reductions.
To implement such mechanisms, we rely on the communication capabilities of the Internet of Things: our mechanisms involve broadcasting a very short message to all appliances, which then need minimal computational effort to respond.

Future work can go in several directions. On the theoretical side, encompassing a variety of appliance types and the possible message losses or delays, as well as exploring more complex coordination schemes, are worth further investigation. On the practical side, the format of the messages to send can be specified, and on-field experiments can be carried out.
Finally, the economic side has not been considered in this paper, but constitutes a major aspect of flexibility markets: our analysis shows how much reduction an aggregator of appliances can offer, but the appliance owners need to be sufficiently incentivized to contribute. Similarly, the flexibility market structure (in particular, the level of competition) will have a strong impact on the prices of reductions and the associated rewards for all participants, and ultimately, on the amounts of flexibility offered by those new means.

\bibliographystyle{unsrt}
\bibliography{biblio}

\end{document}